\documentclass[a4paper,10pt]{amsart}

\usepackage{amssymb,amsthm,amsmath,verbatim}
\usepackage{t1enc}
\usepackage{enumerate}
\usepackage{color}
\usepackage{tikz}

\usepackage{wrapfig}
\usepackage{framed}
\usepackage{lipsum} 
\usepackage[alpine,misc]{ifsym}

\usepackage{mathtools}

\usepackage{xmpmulti}
\usepackage{psfrag}
\usepackage{mathtools}

\usepackage[usenames,dvipsnames]{pstricks}
\usepackage{epsfig}
\usepackage{pst-grad}
\usepackage{pstricks,pst-node}
\usepackage{pspicture}
\usepackage{float}
\usepackage[colorlinks]{hyperref}
\definecolor{mygray}{gray}{0.6}
\usepackage[utf8]{inputenc}
\usepackage[T1]{fontenc}

\usepackage{todonotes}

\newcommand{\mc}[1]{\mathcal{#1}}

\newcommand{\mb}[1]{\mathbb{#1}}

\newcommand{\mf}[1]{\mathfrak{#1}}

\usepackage[top=3cm, left=2.5 cm, right=2.5cm]{geometry}
 \author[D. T. Soukup]{Dániel T. Soukup}
  \address[D.T. Soukup]{Universit\"at Wien,
Kurt G\"odel Research Center for Mathematical Logic, Austria}
 \email[Corresponding author]{daniel.soukup@univie.ac.at}
 \urladdr{http://www.logic.univie.ac.at/$\sim$soukupd73/}

\title{Two infinite quantities and their surprising relationship}

\begin{document}
\maketitle

\begin{wrapfigure}[15]{r}{.5\textwidth}
\centering
 \includegraphics[width=0.4\textwidth]{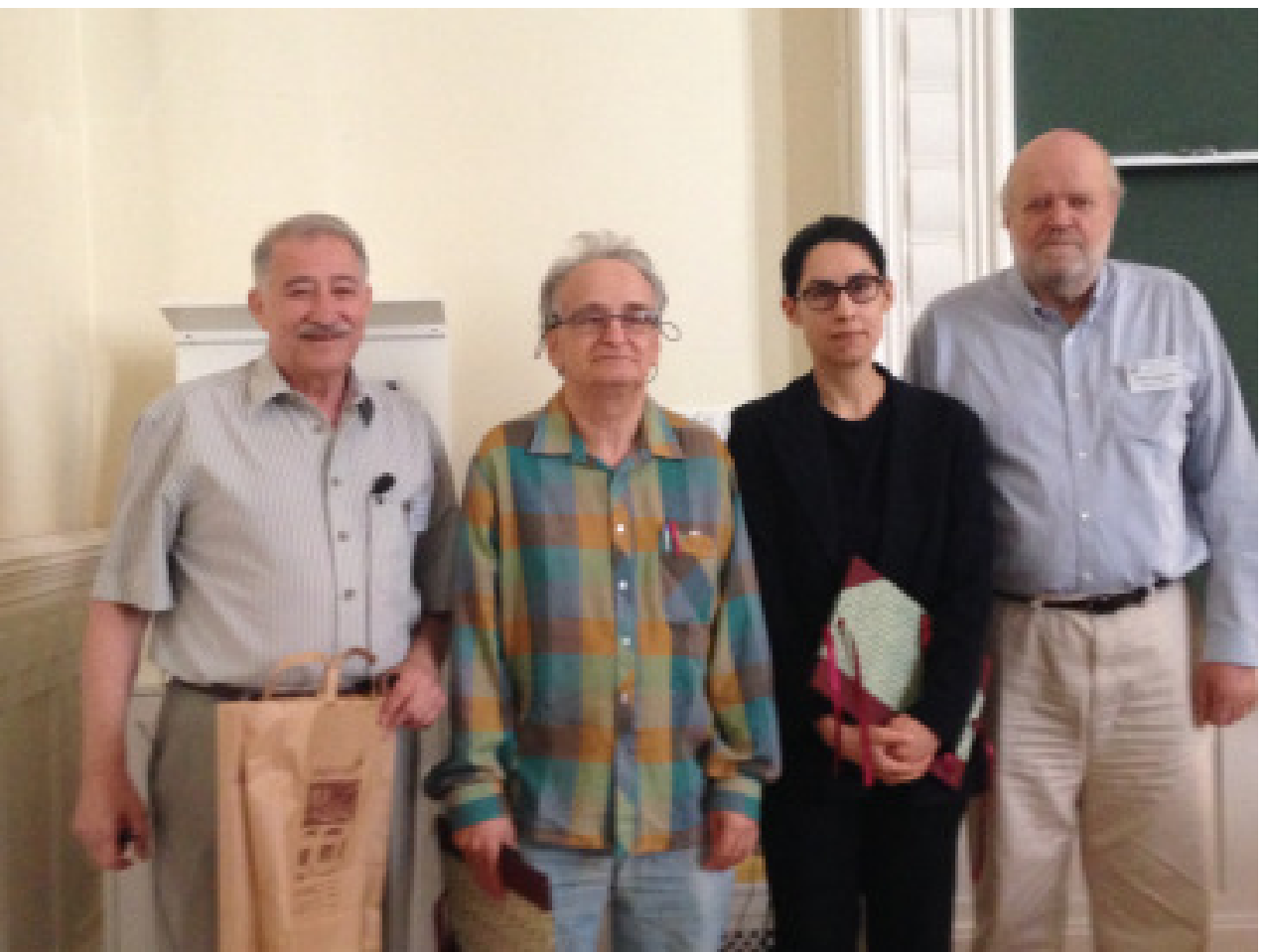}
 
{S. Shelah and M. Malliaris in the center\\

(Photo by Joan Bagaria)}
\end{wrapfigure}


As early as the 17th century, Galileo Galilei wondered how to compare the sizes of infinite sets. Fast forward almost four hundred years, and in the summer of 2017, at the 6th European Set Theory Conference, a young model theorist, Maryanthe Malliaris, and the well-known polymath, Saharon Shelah, received the Hausdorff Medal for the most influential work in set theory published in the last five years. Malliaris and Shelah made significant breakthroughs both regarding a model theoretic classification problem (that is, sorting certain objects into types), and proved that two well-studied infinite quantities, $\mf p$ and $\mf t$, are in fact the same. This latter result is the focus of our expository paper.


\bigskip

Galilei  considered the set of natural numbers  $\mb N$ and the set of perfect squares $\{1,4,9,16\dots\}$. The argument, that these two sets have the same size, goes as follows: since any perfect square has exactly one positive root, and any positive natural number is the root of some perfect square, there should be the same amount of perfect squares and positive natural numbers. On the other hand, there are many natural numbers which are not perfect squares; indeed, if one looks at the ratio of perfect squares to all natural numbers in larger and larger intervals, this quantity tends to zero rather fast. Galilei thought of these observations as a paradox, which prevents us from distinguishing between the sizes of infinite sets.


In the 1870s, Georg Cantor came forward with the following definition, now accepted as standard: \emph{two sets are equinumerous} or \emph{have the same cardinality} exactly if there is a one-to-one correspondence between the elements of the two sets. So, Galilei's first argument proves that  $\mb N$ and the set of perfect squares have the same cardinality. Those sets which are equinumerous with $\mb N$ are called \emph{countably infinite} and we use $\aleph_0$ (in plain words 'aleph zero') to denote their size; the $\aleph_0$   notation refers to the fact that this is the smallest possible infinite size, or cardinality in other words.\footnote{We say that a set $X$ has cardinality \emph{at most the cardinality of $Y$}, if $X$ and some subset of $Y$ are equinumerous.} 


One of Cantor's great contributions to logic was that he did not consider Galileo's argument as an irresolvable paradox, but instead, he started to develop a rich theory of infinities.

\begin{wrapfigure}[11]{r}{.5\textwidth}
\centering


 \includegraphics[width=0.3\textwidth]{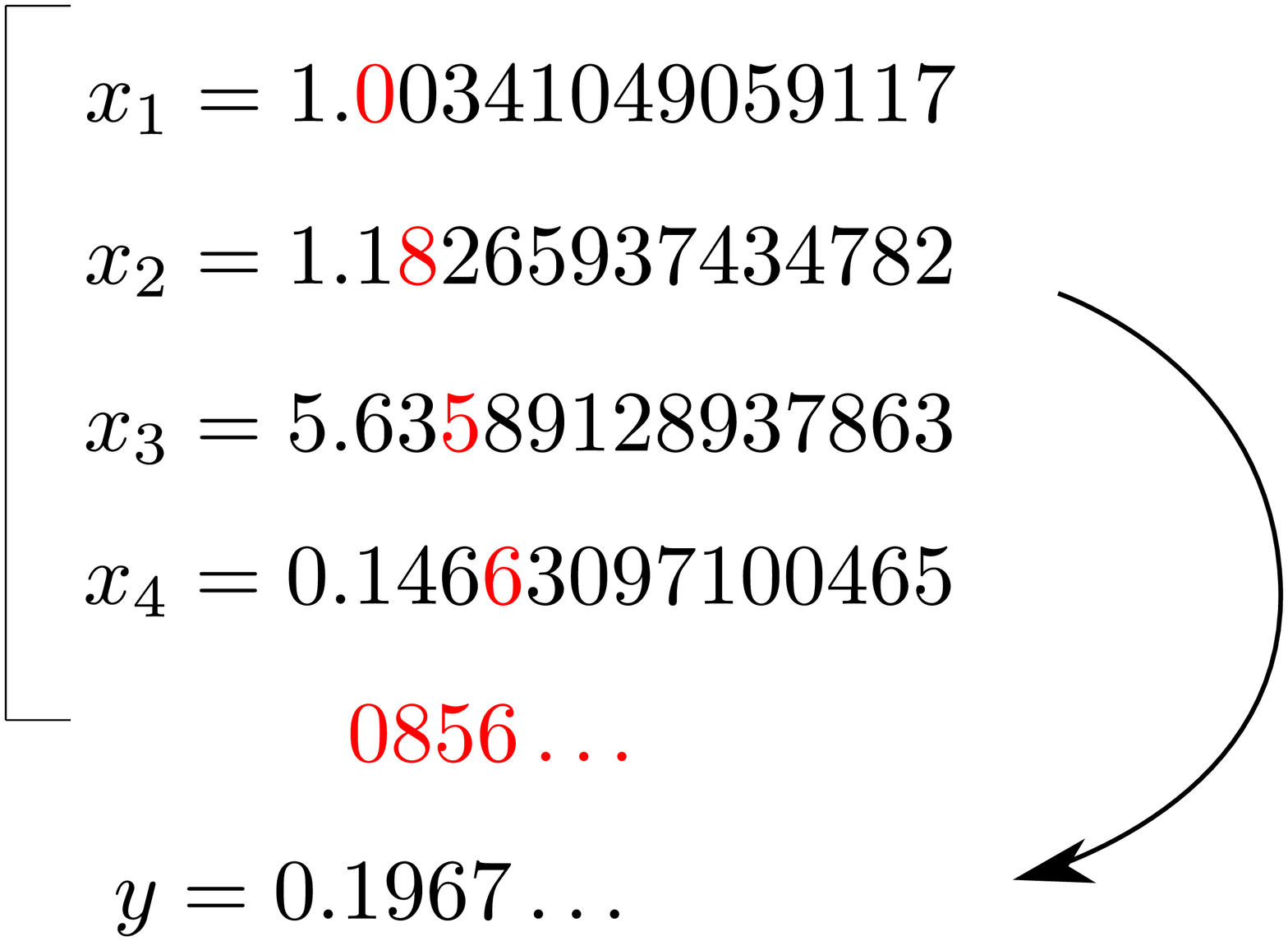}

\end{wrapfigure}

Quite surprisingly, he proved that even the set of rational numbers $\mb Q$ is countably infinite. Now, how about the set of all real number $\mb R$? Cantor claims, that no matter how we produce a list $x_1,x_2,x_3,\dots$ of real numbers, there is always a real number $y$ which is missing from our enumeration.
Indeed, if $y$ differs from $x_1$ at the first decimal place, and $y$ and $x_2$ differ at the second decimal place, and so on, then $y$ cannot possibly appear on the list. Hence, there is no one-to-one correspondence between the natural numbers and the set $\mb R$, and so  $\mb R$ must be \emph{uncountable}. 
We use the notation $2^{\aleph_0}$ (the 'continuum') for the cardinality of $\mb R$, and in turn, we proved the following inequality:
$$ \textmd{the cardinality of }\mb N =\aleph_0<2^{\aleph_0}=\textmd{the cardinality of }\mb R.$$

Based on Cantor's theory, one can compare the cardinalities of any two sets, and in any non empty collection of sets, there is a smallest in size. So, it makes sense to define the first infinite cardinality which is uncountable, and we use $\aleph_1$ to denote this. The next larger cardinality will be denoted by $\aleph_2$, then $\aleph_3$, and so on.\footnote{The list of $\aleph$'s does not stop here: the $\aleph_n$ cardinalities indexed by the natural number have a supremum, denoted by  $\aleph_\omega$, the next strictly larger cardinality is $\aleph_{\omega+1}$, then $\aleph_{\omega+2}$...}





\begin{wrapfigure}[6]{r}{.5\textwidth}
\fbox{\begin{minipage}{\dimexpr\linewidth-3\fboxrule-3\fboxsep}

\textbf{Is there a largest infinite cardinality?} Not according to Cantor's definition. Any set $X$ has more subsets than elements, in other words,\\  $\kappa<2^\kappa$ for any cardinality $\kappa$.


\end{minipage}}
\end{wrapfigure}


\bigskip

At this point, we defined a strictly increasing sequence  $\aleph_0< \aleph_1<\aleph_2<\dots$ of infinite cardinalities. Where does $2^{\aleph_0}$ sit in this list? Interestingly, this question is \emph{undecidable} using the generally accepted ZFC axioms of mathematics.\footnote{ZFC stands for the Zermelo-Fraenkel axiom system with the 'Axiom of Choice'.} Put it differently, in some models of mathematics the equation  $2^{\aleph_0}=\aleph_1$ holds i.e., $\mb R$ has the smallest possible uncountable cardinality, and in this case, we say that the \emph{Continuum Hypothesis} holds. In many other interesting and important models, the Continuum Hypothesis fails and the gap  between ${\aleph_0}$ and $2^{\aleph_0}$ can be arbitrary large: some models satisfy ad hoc equations, such as $2^{\aleph_0}=\aleph_{16}$.


A simple analogy might shed more light on this undecidability business:  $\mb R$ and $\mb Q$, as algebraic structures, both satisfy certain axioms of addition and multiplication,\footnote{Think about commutativity, distributivity or the \textit{field axioms} in general.} however the equation $x^2=2$ has two roots in $\mb R$, but none in $\mb Q$. In other words, whether the statement '$x^2=2$ has a solution' is true or false, is not decided by the axioms of addition and multiplication. Similarly, while the ZFC axiom system decides that 'the inner angles of a plane triangle sum to 180 degrees', the axioms do not decide  whether $2^{\aleph_0}=\aleph_1$ or $2^{\aleph_0}>\aleph_1$ holds; one of these statements holds in any model, but it depends on the particular model which. Another quite similar story was unfolding with the discovery of non-standard geometries in the 19th century and the independence of the \emph{parallel postulate} from other geometric axioms. 
 
 \medskip 

\begin{wrapfigure}[18]{l}{.5\textwidth}
\fbox{\begin{minipage}{\dimexpr\linewidth-3\fboxrule-3\fboxsep}

\textbf{How does one make new models?} Kurt Gödel, in 1938, discovered the so-called \textit{constructible universe}, and proved that the Continuum Hypothesis holds in this model of mathematics. Then, in the 1960s, Paul Cohen genuinely surprised the set theoretic community: he was the first to show that, in some other models of ZFC, the Continuum Hypothesis fails i.e., $\aleph_1<2^{\aleph_0}$. Cohen's technique is called \emph{forcing}, and its underlying idea is fairly simple: given a set theoretic model $M$, we can construct a larger model $N$ by adding a generic object $G$ to $M$. We can use the generic $G$ to increase the value of $2^{\aleph_0}$, which leads to the failure of the Continuum Hypothesis. Cohen received a Fields Medal for his work in set theory in 1966.


\end{minipage}}
\end{wrapfigure}
 




In the last fifty years, people constructed various different models of mathematics, which have plenty of cardinalities between $\aleph_0$ and $2^{\aleph_0}$, and a significant portion of modern set theory focuses on the analysis of these models. Now, two models of the ZFC axioms which both satisfy $2^{\aleph_0}=\aleph_{16}$ can behave very differently, even if one only considers theorems from algebra, measure theory or topology. To the surprise of many specialists, the Whitehead-problem in group theory or the existence of outer automorphisms on the Calkin-algebra are both undecidable using the usual axioms. However, in light of present techniques, it is hard to see how say the Riemann-hypothesis or the existence of general solutions for the Navier-Stokes equations could be undecidable.\footnote{Anyone who solves one of these problems, will be awarded the Millenium Prize and its one million dollar prize purse.}


\medskip

On the bright side, the study of the real line and its fundamental topological and measure theoretic properties is possible through looking at combinatorial properties of families of subsets of the natural numbers; indeed, there is a one-to-one correspondence between the real numbers and subsets of $\mb N$.\footnote{In other words, the cardinality of all subsets of $\mb N$ is $2^{\aleph_0}$.} Studying the so-called \emph{cardinal characteristics of the continuum} deals with exactly such matters.




\bigskip

What attracts many people to this area of research is that one can understand  a breakthrough result such as Malliaris and Shelah's with minimal background. From now on, we will only talk about subsets of the natural numbers, and our concern is the following relation: if $A$ and $B$ are sets of natural numbers, then we write $A\subseteq^* B$ (in plain words, '$A$ is almost contained in $B$') if all but finitely many elements of $A$ are elements of $B$ as well. For example,  the set $A=\{1,3,5,7,\dots \}$ is almost contained in $B=\{4,5,6,7,\dots\}$, since each element of $A$, apart from $1$ and $3$, is an element of $B$ too.


What is the advantage of working with such a weak relation instead of the real containment? Lets take the sets $A_n$ of positive natural numbers which are divisible by $n$: so  $A_1$ is the set of all positive natural numbers, $A_2$ collects the even natural numbers, $A_3=\{3,6,9\dots\}$, and so on.



\begin{wrapfigure}[12]{l}{.45\textwidth}
\centering

 \includegraphics[height=4
 cm]{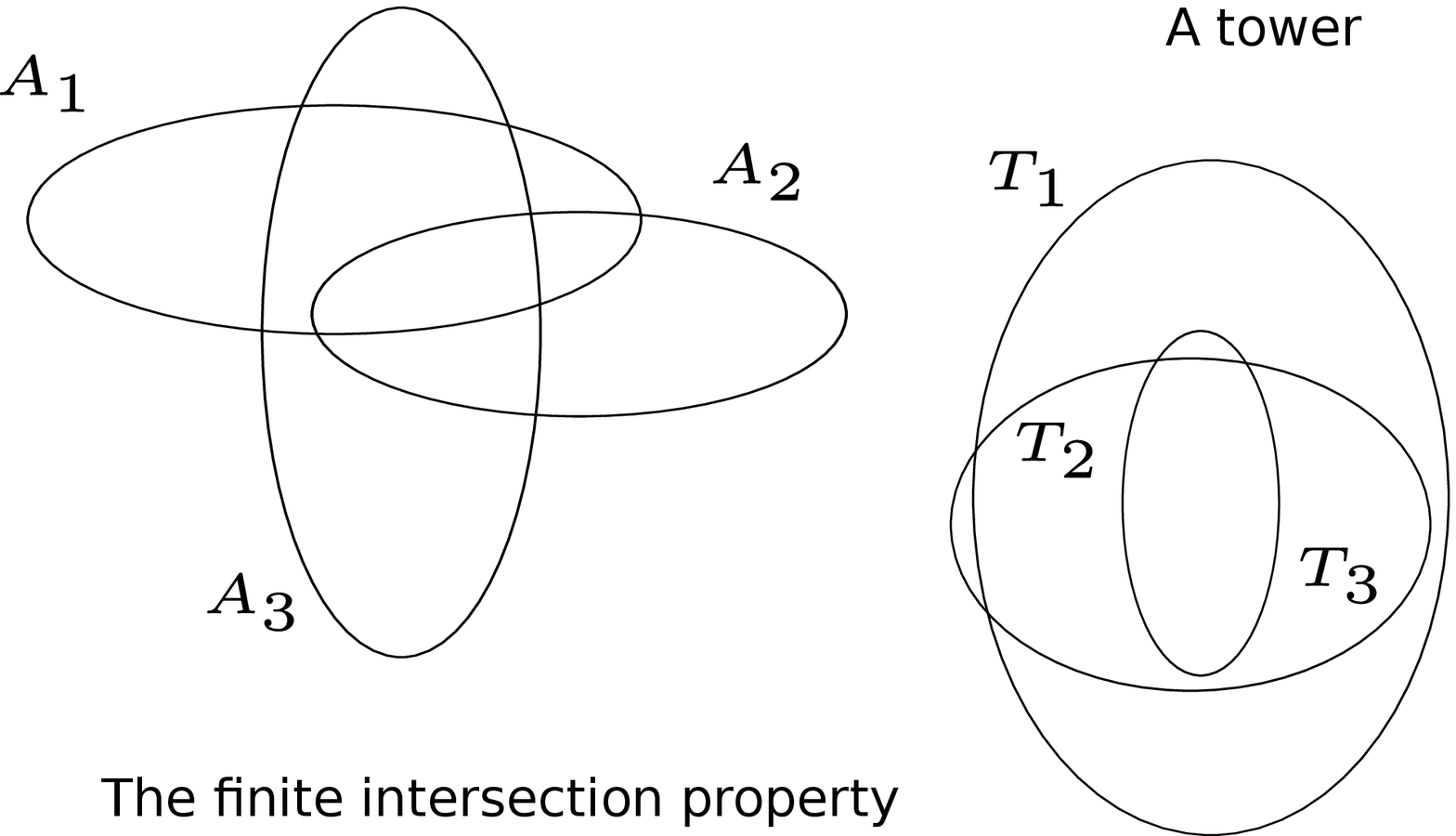}
 \vspace{0.5cm}
 
\end{wrapfigure} 
It is easy to see that if we take finitely many sets  $A_1, A_2\dots A_n$ until a fixed $n$, then these sets have infinite intersection. Indeed, the multiples of $n!=1\cdot 2\cdot\dots \empty\cdot n$ are contained in each of the sets $A_1, A_2\dots A_n$.
 We usually say in this case that the system of sets $\{A_n\}$ has the \emph{finite intersection property}, a rather counterintuitive name for sets with infinite intersection...
 

 Of course we cannot find any positive natural number which is divisible by all the numbers, that is, there is no real intersection to the whole family $\{A_n\}$. On the other hand, we can easily find infinite sets $B$ so that $B\subseteq^* A_n$ for all $n$; it suffices to make sure that the $n$th element of $B$ is selected from  $A_1\cap A_2\cap \dots \cap A_n$, so for example $B=\{n!\}_{n=1,2\dots}$ works. Since $B$ behaves somewhat like an intersection to the family $\{A_n\}$, we call $B$ a \emph{pseudo-intersection} of $\{A_n\}$. It is also easily checked, that one can add $B$ to the family 
 $\{A_n\}$, and the finite intersection property is still preserved.

 
 
 Now, given an arbitrary family of sets $\mc A$ with the finite intersection property, we can extend $\mc A$ to a maximal family $\mc A_{\max}$ with the finite intersection property.\footnote{This is a standard application of Zorn's lemma.} This maximal family $\mc A_{\max}$ however cannot have a pseudo-intersection anymore.\footnote{Otherwise, we could add this pseudo-intersection to $\mc A_{\max}$ while preserving the finite intersection property, and so maximality would be violated.} This leads to our first main definition: 
\bigskip
 \begin{center}
  
 \begin{minipage}{0.5\textwidth}
   The \emph{pseudo-intersection number} $\mf p$ is the cardinality of the smallest family of sets $\mc A$ with the finite intersection property, such that $\mc A$ does not have a pseudo-intersection. 
 
 \end{minipage}

 \end{center}
\bigskip

 
 
 In the above example, we essentially showed that any family \emph{indexed by the natural numbers}, which has the finite intersection property, also has a pseudo-intersection, and in turn, $\mf p$ must be uncountable. Hence, the following inequality holds: $$\aleph_0<\mf p\leq 2^{\aleph_0}.$$
 
 
 
 There are numerous models of mathematics, where  $\mf p= 2^{\aleph_0}$ and this common value can be essentially any $\aleph$. On the other hand,  $\aleph_1=\mf p<2^{\aleph_0}=\aleph_2$ can also be true in other models. So the usual axioms do not decide where the pseudo-intersection number  $\mf p$ sits on the list  $\aleph_1,\aleph_2,\dots$, nor if $\mf p= 2^{\aleph_0}$ or $\mf p< 2^{\aleph_0}$ holds.
 
 

\bigskip


 We need another definition. A typical family with the finite intersection property is not ordered in any sense: in our original example about divisibility, if one only considers the sets $A_p$ for primes $p$, no two of these sets are in $\subseteq^*$ relation. Now, call a family of sets $\mc T$ a \emph{tower} if for any two sets $X,Y$ from $\mc T$, either $X\subseteq^* Y$ or $Y\subseteq^* X$ holds. In other words, the relation $\subseteq^*$ linearly orders  $\mc T$. 
 
 \bigskip
  \begin{center}

 \begin{minipage}{0.5\textwidth}
  
 The so called \emph{tower number} $\mf t$ is the cardinality of the smallest tower  $\mc T$ without pseudo-intersection.
  \end{minipage}

 \end{center}
\bigskip


\begin{wrapfigure}[14]{r}{.5\textwidth}
\fbox{\begin{minipage}{0.5\textwidth}

\textbf{Large towers.} Surprisingly, even from the subsets of $\mb N$, one can construct towers of size  $2^{\aleph_0}$. First of all, list all the rational numbers as $q_1,q_2,q_3,\dots $. Then, for each real number $r$, we define $X_r$ to be the set of all those indexes $n$ so that $q_n<r$. So $X_r$  essentially collects the rational numbers smaller than $r$.  Now, if $r<t$ are two real numbers then  $X_r$ is fully contained in $X_t$, moreover, 
$X_t$ has infinitely many extra elements. Can we find a pseudo-intersection for this family?


\end{minipage}}
\end{wrapfigure}

Since any tower has the finite intersection property, any witness for the invariant $\mf p$ must have at most the size of a tower witnessing $\mf t$. This shows the next equation: $$\hspace{-5cm} \aleph_0<\mf p\leq \mf t\leq 2^{\aleph_0}.$$



%

 \smallskip

 The value of the invariant $\mf t$ can be manipulated similarly to $\mf p$. Actually, since the 1940s, more than a dozen cardinal invariants between $\aleph_0$ and $2^{\aleph_0}$ have been studied. Moreover, we know that apart from certain simple inequalities that were known from the mid 20th century, there is no provable relationship between the invariants. That is, with various versions and combinations of forcing techniques, one can not only set the value of $2^{\aleph_0}$ but also the values of the invariants to any fixed alephs, which do not violate the known inequalities.\footnote{We recommend A. Blass' classical text on cardinal invariants for an excellent overview.}


\begin{wrapfigure}[17]{l}{.45\textwidth}
\centering

   \vspace{-0.3cm}
 \includegraphics[height=5.5cm]{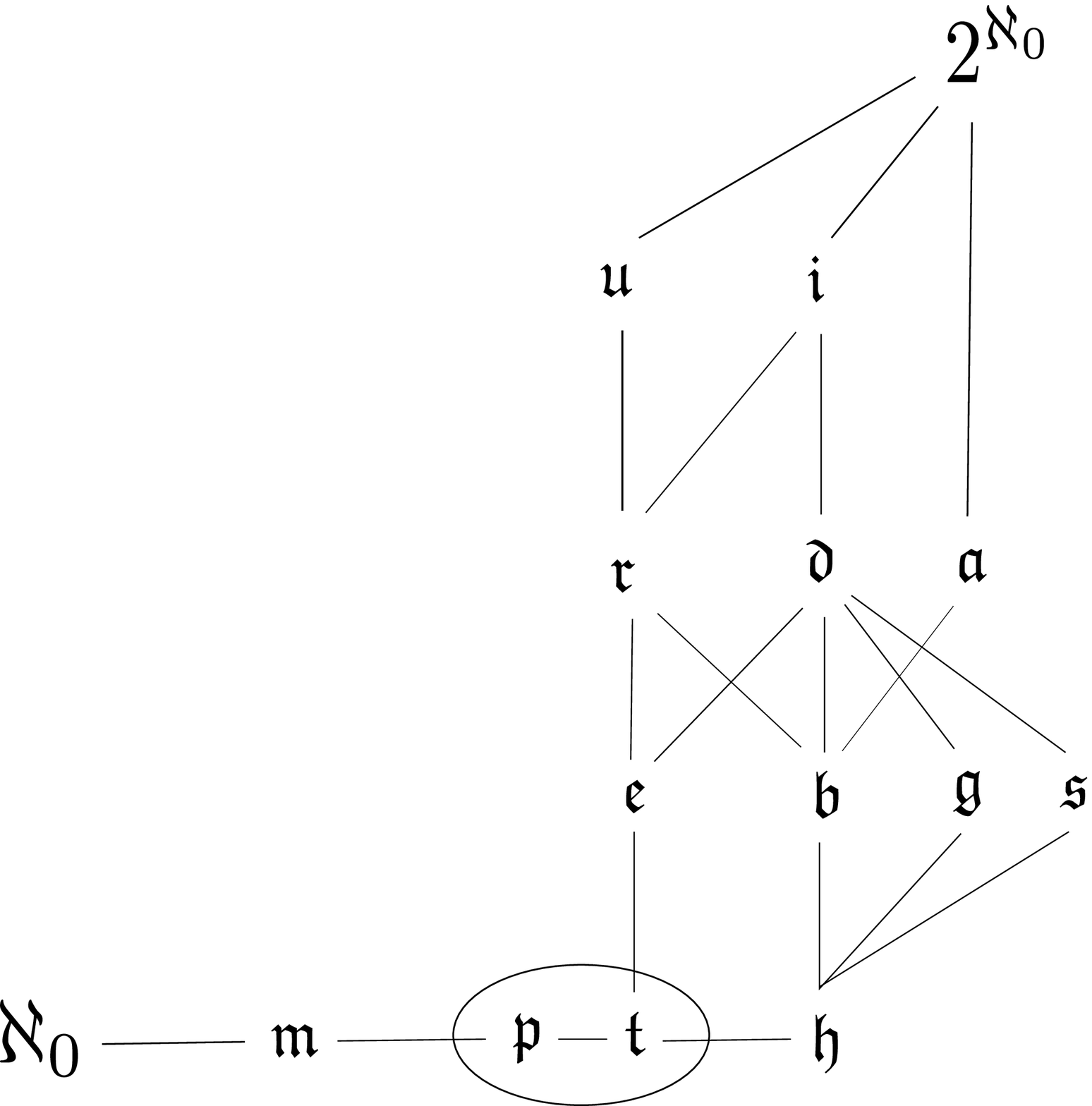}
  \vspace{0.2cm}
  
  {A few cardinal invariants and their provable relationships}

\end{wrapfigure}

Indeed, the most advanced techniques can independently manipulate the values of more than five invariants at the same time.


\bigskip 

Despite all these investigations in the last sixty years, no one constructed a model, where $\mf p$ and $\mf t$ are not equal. The generally accepted conjecture was that $\mf p<\mf t$ is possible, however, we knew that a proof could not be too easy: it was proved a long time ago that if   $\mf p=\aleph_1$, then $\mf t=\aleph_1$ holds as well. So, $\mf p<\mf t$ can only hold in models where $2^{\aleph_0}$ is bigger than $\aleph_2$, however fine tuning such models is significantly harder than controlling models with $2^{\aleph_0}\leq \aleph_2$.

 
 \medskip

 Malliaris and Shelah's new, unexpected result states that $$\hspace{8cm }\mf p=\mf t,$$ no matter which model we look at. What is it exactly that the authors prove? Since $\mf p\leq \mf t$ was known, they needed to show $\mf t\leq \mf p$: given an arbitrary family $\mc A$ with the finite intersection property but no pseudo-intersection, no matter how random or involved the overlays between the elements of $\mc A$ are, one can construct a tower  $\mc T$ of size at most the size of $\mc A$ which has no pseudo-intersection. The naive approach to construct the tower from the elements of $\mc A$ fails quickly, since  $\mc A$ might have no two elements which are related by $\subseteq^*$, but $\mc T$ must be totally ordered by $\subseteq^*$.


One can compare Malliaris and Shelah's result to showing that two algebraic equations have the same solutions, but without actually determining the value of these solutions. The ZFC axioms do not decide whether $\mf p=\aleph_1$ or  $\mf p=\aleph_2$ or  $\mf p=\aleph_3$, and similarly the value of $\mf t$ is undecidable. However, Malliaris and Shelah could prove that no matter what the value of $\mf p$ is, it has to be the same as $\mf t$.  



\bigskip



 We should mention that not only the authors solved the sixty year old mystery surrounding $\mf p$ and $\mf t$, but they uncovered a novel connection between a model theoretic complexity hierarchy, the Keisler-order, and the theory of cardinal characteristics \cite{MS1,MS2, moorept}; unfortunately, it is not in the scope of our paper to sketch these results. While the original proof of $\mf p=\mf t$ employs serious tools from model theory and modern set theory, there are now new versions which only require basic knowledge and some perseverance
 \cite{ptthesis}\footnote{The Fields Medalist, Timothy Gowers' blog also looks at this problem, \href{https://gowers.wordpress.com/2017/09/19/two-infinities-that-are-surprisingly-equal/}{link here}.}

\bigskip

\begin{wrapfigure}[16]{l}{.5\textwidth}
\fbox{\begin{minipage}{0.45\textwidth}

\textbf{An open problem.} We say that a family  $\mc R$ is \emph{unsplit} if there is no $Y$ which splits all elements of $\mc R$ at the same time, i.e., no $Y$  so that the intersection $X\cap Y$ and the difference $X- Y$ are both infinite for any $X$ from $\mc R$. Let  $\mf r$ denote the cardinality of the smallest unsplit family. Moreover, let $\mf r_{\sigma}$ denote the size of the smallest family which cannot be split by even countably many sets  $Y_0,Y_1,\dots$. It is easy to see that $\aleph_0< \mf r\leq \mf r_{\sigma}\leq 2^{\aleph_0}$, however, it remains unknown if $\mf r< \mf r_{\sigma}$ is possible in some model. The conjecture is that $\mf r= \mf r_{\sigma}$ holds, and this is certainly true in all known models \cite{splitting}.

\end{minipage}}
\end{wrapfigure}


The results of Malliaris and Shelah are far from the last of cardinal characteristics, and will more likely spark a renewed interest in the field. So what problems does a regular set theorist work on? On one hand, the relationship of some classical invariants are still unknown, and we mention a problem of this sort in the side note \cite{splitting}. On the other hand, people are  defining new, interesting invariants to this day, and it is often a hard task to determine the position of these new invariants relative to the classical ones \cite{rearrangement}. Finally, a rich theory is growing out of the study of cardinal characteristics which are defined using families of uncountable sets, rather than the subsets of  $\mb N$, showing striking differences with the classical studies \cite{brendle}.

 
\bigskip

\bigskip

\vspace{-0.3cm}

We close by a few words about the awardees: Maryanthe Malliaris graduated from Berkeley in 2009, currently  a professor at the University of Chicago,  she is the recipient of multiple, prestigious awards, and is invited to present at the  International Congress of Mathematicians in 2018.

\medskip
 
The name Saharon Shelah might ring a bell for a lot of the readers: the 72 years old mathematician is the author of 1023 published papers (!), on groundbreaking results from combinatorics and model theory, to logic and group theory. He still works 6 days a week, splitting the year between Rutgers and the Hebrew University of Jerusalem.
 





\vfill

{\Small
\emph{The author was supported in part by the FWF Grant I1921. A Hungarian version of the current survey was prepared for the journal \emph{Matematikai Lapok}. We thank  Emese Bottyán, Lajos Soukup, Zoltán Vidnyánszky, and Zita Zádorvölgyi for their careful reading, while the English version significantly improved thanks to the help of Neil Barton and Vera Fischer.}
}




\end{document}